\documentclass[12pt]{article}
\usepackage[all]{xy}
\usepackage{amsmath}
\usepackage{amsfonts}
\usepackage{amssymb}
\usepackage{amscd}
\usepackage{amsthm}
\usepackage{latexsym}
\usepackage{amsbsy}

\newtheorem{lem}{Lemma}[section]
\newtheorem{prop}{Proposition}[section]
\newtheorem{theorem}{Theorem}[section]
\newtheorem{corollary}{Corollary}[section]
\newtheorem{definition}{Definition}[section]
\newtheorem{example}{Example}[section]

\newcommand{\ot}{\otimes} 
\newcommand{\vs}[2]{\vartheta_{#1}^{(#2)}}
\begin{document}
\title
{A New Cyclic Module for Hopf  Algebras \footnote{The results of this paper were announced in the AMS meeting in
 Toronto, September 23-24, 2000. }}

\author {M. Khalkhali,~~~ B. Rangipour,
\\\texttt{~masoud@uwo.ca ~~~~brangipo@uwo.ca}
\\ Department of Mathematics 
\\ University of Western Ontario }
\date{}
\maketitle

\begin{abstract}
We define a new cyclic module, dual to the Connes-Moscovici cocyclic module, for Hopf algebras, and give a
 characteristic  map for  coactions of Hopf algebras. We also compute the resulting cyclic homology for cocommutative
 Hopf algebras, and some quantum groups.
 \end{abstract}


\section{Introduction}
In their study of the index theory of transversally elliptic operators $~\cite{aChM98}$,  Connes and Moscovici
 developed a cyclic (co-)homology theory for Hopf algebras, which can be considered as an extension  of group homology
  and Lie algebra homology to Hopf algebras and in particular  to quantum groups. This theory was further explained
   and developed along purely algebraic lines in ~$\cite{aChM99,cr}$.  One of the main tools in
    $~\cite{aChM98,aChM99}$
 is a noncommutative characteristic map ${HC}_{(\delta,\sigma)}^\ast (\mathcal{H}) \longrightarrow  HC^\ast (A)$,
  for a Hopf
  algebra $\mathcal{H}$ and a $\mathcal{H}$-module algebra $A$ endowed with an invariant trace.

There is however a need for a dual theory to be developed. This is justified for example when one studies coactions of
 Hopf algebras and also by the fact, first observed by M. Crainic, that for group algebras and 
  in general for Hopf algebras 
 with a normalized Haar integral (i.e. cosemisimple Hopf algebras),
   that both cyclic homology and cohomology are trivial~$\cite{cr}$.

In this paper we define a new  cyclic module for Hopf algebras. We define the characteristic map for the coaction of
 Hopf algebras and also prove an analogue of Karoubi's theorem for cocommutative Hopf algebras. We also compute our
  theory for quantum groups $U_q(sl_2)$ and $A(SL_q(2))$. It would be very interesting to compute this theory for
   other quantum groups.
 In \cite{cr} one can find a computation of Hopf algebra cyclic cohomology,  in the sense of Connes-Moscovici, for the
  quantum group $U_q(sl_2)$. The method however is very different from ours. 
   We note that a similar cyclic module is also independently considered by R. Taillefer $~\cite{taill}$, for
    different reasons.

    We would like to thank the referee whose suggestions and comments improved our presentation,
     specially in the last section where our original formulation of  Theorem 5.1 and Corollary 5.2 were inaccurate.

\section{Cyclic Module of Hopf Algebras}
Let $(\mathcal{H},\mu ,\eta ,\delta  , \epsilon ,S )$ be a Hopf algebra over a commutative 
ring $k$, where  $\mu ,\eta
 ,\delta  , \epsilon ,S $ denote the product, unit map, coproduct, counit and the antipode
 of $\mathcal{H}$,
  respectively. Let  $\sigma$ be a grouplike element of $\mathcal{H}$
and $\delta  : \mathcal{H} \longrightarrow   k $  be a character for $\mathcal{H}$ as
 in$~\cite{aChM99}$. The pair $(\delta  , \sigma ) $ is called 
a modular pair if $\delta ( \sigma) =id $, and a modular pair in involution if 
\begin{equation} 
\widetilde {S}_{\sigma}^2=id,
\end{equation}
where \[\widetilde {S}_{\sigma}(h)=\sigma \sum_{(h)}\delta (h^{(2)})S(h^{(1)}). \] 
We have used  Sweedler's notation$~\cite{sw}$
i.e., $ \Delta (h)= \sum_{(h)}h^{(1)} \otimes h^{(2)}. $
We will associate a cyclic module to any Hopf algebra $\mathcal{H}$ over $k$ if $\mathcal{H}$ has a modular pair
 $(\delta  , \sigma ) $
in involution. This cyclic module somehow can be seen as the dual of the cocyclic module  introduced in
 $~\cite{aChM99,aChM98}$ by A. Connes and H. Moscovici. First, consider $S_{\sigma}$  (where $ S_{\sigma} (h)= \sigma
  S(h) $ )  which has the properties:
$$ S_{\sigma} (h_1h_2)=S_{\sigma}(h_2)S(h_1) $$
$$ S_{\sigma}(1)=\sigma $$
\[ \Delta S_{\sigma}(h)= \sum_{(h)} S_{\sigma}(h^{(2)})\otimes S_{\sigma}(h^{(1)}) \] 
$$ \epsilon (S_{\sigma}(h))=\epsilon(h). $$

Using $\epsilon$ and $\delta  $ one can endow $k$ with  an $\mathcal{H}$-bimodule structure, i.e.,
$$\delta  \otimes id : \mathcal{H} \otimes k \longrightarrow  k \quad and \quad 
 id \otimes \epsilon : k \otimes \mathcal{H} \longrightarrow  k.
 $$
Our cyclic module as a simplicial module is exactly the 
Hochschild complex of the algebra $ \mathcal{H} $ with coefficients in $ k $,  where $ k $ is an 
 $\mathcal{H}$-bimodule as  above. So  if  we  denote our cyclic module by
${\widetilde {\mathcal{H}}^{(\delta,\sigma)}}, $  we have
 $ {\widetilde {\mathcal{H}}}^{(\delta,\sigma)} _n = \mathcal{H}^{\otimes n}$and 
   $ {\widetilde {\mathcal{H}}}^{(\delta,\sigma)} _0 =k$.
  Its 
faces and degeneracies are as follows: \\ 
\begin{eqnarray*}
 {\delta }_0 (h_1 \otimes h_2 \otimes  ... \otimes h_n) &=&\epsilon (h_1)h_2 \otimes h_3 \otimes ... \otimes h_n \\
 {\delta }_i (h_1 \otimes h_2 \otimes 
 ... \otimes h_n) &=& h_1 \otimes h_2 \otimes ...\otimes 
h_i h_{i+1}\otimes ... \otimes h_n~~  1\le i\le n-1 \\
{\delta }_n(h_1 \otimes h_2 \otimes  ... \otimes h_n) &=&\delta  (h_n)h_1\otimes h_2 \otimes ... \otimes h_{n-1} \\
 {\sigma}_0(h_1 \otimes h_2 \otimes  ... \otimes h_n) &=&1 \otimes h_1  \otimes ... \otimes h_n \\
\hspace{2cm}{\sigma}_i(h_1 \otimes h_2 \otimes ... \otimes h_n) &=& h_1 \otimes h_2 ...\otimes h_i \otimes 1 \otimes
 h_{i+1} ... \otimes h_n  ~~1\le i\le n-1 \\
 {\sigma }_n(h_1 \otimes h_2 \otimes  ... \otimes h_n) &=& h_1 \otimes h_2  \otimes ...\otimes    h_n\otimes 1. 
\end{eqnarray*} 
To define a cyclic module it remains to introduce an action of the cyclic group on our module. Our candidate is 
$$\tau _n(  h_1 \otimes h_2 \otimes  ... \otimes h_n ) = \sum\delta  (h_n^{(2)})
( S_{\sigma} ( h_1^{(1)} h_2^{(1)} ...h_n^{(1)} ) \otimes h_1^{(2)} \otimes ...\otimes h_{n-1}^{(2)} ). $$
\begin {theorem}
Let $(\mathcal{H},\mu ,\eta ,\delta  , \epsilon , S ) $ be a Hopf algebra 
over $k$  with a modular pair $(\delta  ,
 \sigma)$ in involution. 
Then $\widetilde{\mathcal{H}}^{(\delta,\sigma)}$
 with operators given  above  defines a cyclic module. Conversely, 
 if $\delta  (\sigma ) = 1$ and
  $\widetilde{\mathcal{H}}^{(\delta,\sigma)}$ is a cyclic module,
   then $(\delta  , \sigma)$ is a modular pair in involution. 
\end {theorem}
\begin{proof}  As we mentioned the simplicial relations are already held and it remains to check 
the following extra relations:
\begin{eqnarray}
{\tau}_n^{n+1} &=& id \\
{\delta }_i {\tau}_n &=& {\tau}_{n-1}{\delta }_{i-1} \\
\delta _0 {\tau}_n & =& {\delta }_n \\
{\sigma}_i {\tau}_n &=& {\tau}_{n+1}{\sigma}_{i-1} \\
{\sigma}_0 {\tau}_n &=& {\tau}_{n+1}^2{\sigma}_{n}. 
\end{eqnarray}

To prove equation (2), we first compute $\tau _n^2$ :
\begin{eqnarray*}
 \tau_n^2(h_1 \otimes h_2 \otimes ... \otimes h_n ) &=& \tau ( \sum\delta (h_n^2)S_\sigma (h_1^{(1)}h_2^{(1)}...
  h_n^{(1)})\otimes h_1^{(2)}\otimes ...\otimes h_{n-1}^{(2)} ) \\
 &=&\sum\delta  (h_n^{(2)})\tau_n(S_\sigma (h_1^{(1)}h_2^{(1)}... h_n^{(1)})\otimes h_1^{(2)}\otimes ...\otimes
  h_{n-1}^{(2)} ) \\
 &=& \sum\delta  (h_n^{(2)})\sum\delta 
  ((h_{n-1}^{(2)})^{(2)})S_{\sigma}(S_\sigma(h_1^{(1)}...h_n^{(1)})^{(1)}(h_1^{(2)})^{(1)}...\\&&
(h_{n-1}^{(2)})^{(1)})\otimes  
 (S_\sigma(h_1^{(1)}...h_n^{(1)})^{(2)}\otimes (h_1^{(2)}))^{(2)}\otimes ...\otimes (h_{n-2}^{(2)})^{(2)} \\
 &=& \sum\delta  (h_n^{(3)})\delta  (h_{n-1}^{(4)})S_{\sigma}( \sigma
  S(h_n^{(2)})...S(h_1^{(2)})h_1^{(3)}...h_{n-1}^{(3)}) \otimes \\&& 
 S_{\sigma}(h_1^{(1)}h_2^{(1)}...h_n^{(1)})\otimes h_1^{(4)} ...\otimes h_{n-2}^{(4)} \\
&=&\sum\delta (h_n^{(3)})\delta  (h_{n-1}^{(3)})S_{\sigma}( \sigma S(h_n^{(2)})\epsilon (h_{n-1}^{(2)})...\\
&&\epsilon(h_1^{(2)})) \otimes S_{\sigma}(h_1^{(1)}...h_{n}^{(1)}) 
\otimes  
 (h_1^{(3)}\otimes h_2^{(3)}\otimes ...\otimes h_{n-2}^{(3)}) \\
 &=&\sum\delta (h_{n}^{(3)})\delta  (h_{n-1}^{(2)})S_{\sigma }^2(h_n^{(2)})\otimes S_{\sigma
  }(h_1^{(1)}...h_{n}^{(1)}) \\&&\otimes h_1^{(2)}\otimes h_2^{(2)}\otimes ...\otimes h_{n-2}^{(2)}. 
\end{eqnarray*}
 By a similar argument, we can deduce
\begin{multline*}\tau_n^3(h_1 \otimes h_2 \otimes ...\otimes h_n ) = \sum\delta  (h_n^{(3)})\delta 
 (h_{n-1}^{(3)})\delta  (h_{n-2}^{(2)}) 
\\ S_{\sigma}^2(h_{n-1}^{(2)}) \otimes S_{\sigma}^2(h_{n}^{(2)})\otimes S_{\sigma}(h_1^{(1)}...h_n^{(1)}) \otimes
 h_1^{(2)}\otimes h_2^{(2)}\otimes ...\otimes h_{n-3}^{(2)}.
\end{multline*} 
Continuing,
\begin{multline*} \tau_n^{n}(h_1 \otimes h_2 \otimes ...\otimes h_n ) \\= \sum\delta  (h_n^{(3)}) ...\delta
 (h_{2}^{(3)})\delta  (h_{1}^{(2)}) S_{\sigma}^2(h_2^{(2)})\otimes 
                             ... S_{\sigma}^2(h_n^{(2)})\otimes S_{\sigma}(h_1^{(1)}...h_n^{(1)})
\end{multline*}  
and eventually,
\begin{multline*} \tau_n^{n+1}(h_1 \otimes h_2 \otimes ...\otimes h_n) = \\ \sum\delta  (h_n^{(3)}) ...\delta
 (h_{1}^{(3)})\delta  (S(h_{1}^{(1)}))...\delta (S(h_n^{(1)}))S_{\sigma}^2(h_1^{(2)})\otimes ...\otimes
 S_{\sigma}^2(h_n^{(2)}) 
\end{multline*}
$$ = \widetilde {S}_\sigma^2(h_1)\otimes ...\otimes \widetilde {S}_\sigma^2(h_n) 
 = h_1 \otimes h_2 \otimes ...\otimes h_n.$$                        
We leave it to the reader to check the remaining equations. To prove the converse it suffices to 
use just $\tau_1^2=id $.
\end{proof}
\begin{example}
Let ${G}$ be a discrete group and $k{G}$ be its group algebra over $k.$ It is  
a cocommutative Hopf algebra with the
 coproduct, counit and antipode defined by $\Delta(g)=g \otimes g $; 
  $\epsilon(g)=1$; $S(g)=g^{-1}$. We compute the
 cyclic module $\widetilde{kG}^{(\epsilon,1)}$.  
It is obvious that $ k{G}^{\otimes n}$ can be identified with $ k{G}^n, 
$ the free $k$ module generated by $ {G}^n. $
 So we have 
\begin{eqnarray*}
\delta  _i(g_1 ,...,g_n)= \left \{ \begin{array}{ll}
(g_2 ,...,g_n) & \textrm {if $ i=0 $} \\
(g_1 ,...,g_i g_{i+1},...g_n)\ & \textrm{$1 \leq i < n$}\\
(g_1 ,...,{g}_{n-1}) & \textrm{if $ i=n $}
\end{array}\right.
\end{eqnarray*}
\begin{eqnarray*}
\hspace{0.5cm}\sigma _i(g_1 ,...,g_n)=\left \{ \begin{array}{ll}
(1,g_1 ,...,g_n) & \textrm{if $ i=0 $} \\
(g_1 ,...,g_i,1,g_{i+1},...g_n)\ & \textrm{$1 \leq i \leq n-1$}\\
(g_1 ,...,g_{n-1},g_n,1) & \textrm{if $ i=n $}
\end{array}\right.
 \end{eqnarray*}
$$\hspace{1.7cm}\tau(g_1,g_2,...,g_n)= ((g_1g_2...g_n)^{-1},g_1,...,g_{n-1}).$$
\end{example}

It follows that the cyclic module $\widetilde{k{G}}^{(\epsilon,1)}$ 
exactly coincides with $kB{G},$ the cyclic module
 associated with the classifying space of 
${G}~\cite{ld}$.  
 \begin{corollary}
Let ${G}$ be as in the previous example. Then we have 
\begin{eqnarray*}
\hspace{1cm}\widetilde{HP}^{(\epsilon,1)}_n(k{G})=\lim _{\leftarrow}HC_{n+2i}({kBG})= 
\left \{ \begin{array}{ll}
{\prod _{i\geq 0}H_{2i}({G} ;k)} & \textrm{$n$ even}  \\
{\prod _{i\geq 0}H_{2i+1}({G} ;k)} & \textrm{$n$ odd.}
\end{array}\right.
\end{eqnarray*}
\end{corollary}
\section{Hopf Algebra Coaction on an Algebra}
{\text In this section we consider a Hopf algebra $\mathcal{H}$ that has a right coaction on an algebra $A$. In
 technical terms $A$ is 
a right comodule algebra   i.e., there is a $k$-linear map $\beta :A \longrightarrow  A \otimes \mathcal{H}$ such that 
the following diagrams commute 
$$
\xymatrix{
A \ar[r]^\beta \ar[d]_\beta & A \otimes \mathcal{H} \ar[d]^{I_A \otimes \Delta} \\
A \otimes \mathcal{H}\ar[r]^{\beta\otimes I_\mathcal{H}} & A \otimes \mathcal{H}\otimes \mathcal{H}}
\hspace{30pt}
\xymatrix{
& A \ar[r]^{\beta}\ar[dr]_{\cong} & A \otimes \mathcal{H}\ar[d]^{I_A \otimes \epsilon} \\
&  &A\otimes k  } $$
and  $\beta $ is a algebra map where the algebra structure of $A \otimes \mathcal{H} $ is the tensor product of the
 algebras $A$ and $\mathcal{H}$. Similarly one can define a left comodule algebra.
In a  practical way, we can state the commutativity by 
\begin{eqnarray}
 \sum \sum a_{(0)} \otimes (a_{(1)})^{(1)}\otimes (a_{(1)})^{(2)} &=& \sum \sum (a_{(0)})_{(0)}\otimes
  (a_{(0)})_{(1)}\otimes a_{(1)} \\             
 a &=& \sum a_{(0)}\epsilon(a_{(1)}) 
\end{eqnarray}
where the above notations mean
$$ \beta (a) = \sum a_{(0)}\otimes a_{(1)} $$
$$ \Delta (h) = \sum h^{(1)}\otimes h^{(2)}.$$
\begin{definition} A linear map, $Tr :A\rightarrow k$ is called $\delta$-trace  if 
\begin{equation*}
Tr(ab)= \sum _{(b)} Tr(ba^{(0)})\delta  (a^{(1)}) \hspace{2cm}\forall a , b\in  A.
\end{equation*} 
It is called $\sigma $-invariant if for  all $a , b\in  A$,
\begin{multline*}
 ~~~~~~~~~~~~~~~~~~
  \sum _{(b)} Tr(a^{(0)}b)\ (a^{(1)})=\sum _{(a)}Tr(ab^{(0)})S_{\sigma }(b^{(1)}),\\
    \shoveleft{\text{or equivalently }}\\
    \shoveleft{~~~~~~~~~~~~~~~~~~~~ Tr(a^{(0)})a^{(1)}=Tr(a)\sigma.}\\
\end{multline*}
\end{definition}
Let $C_\ast(A)$ denote the cyclic module of the algebra $A.$
\begin{prop}
Let $Tr$ be a $\delta $-trace on $A$, which is  $\sigma$-invariant.  
Then the following map is a cyclic map:
$$ \gamma : C_n(A)  \longrightarrow   \widetilde { \mathcal{H}}_n^{(\delta,\sigma)}, $$
\begin{equation}
 \gamma (a_0 \otimes a_1 \otimes ... \otimes a_n) = \sum Tr(a_0 a_1^{(0)}a_2^{(0)} ... a_n^{(0)})(a_1^{(1)}\otimes... 
  \otimes a_n^{(1)}). 
\end{equation}
\end{prop}
\begin{proof}

We should verify that $\gamma$ commutes with $\delta  _i, \sigma_i, \tau_n:  $ \\
\begin{eqnarray*}
\gamma \circ\delta  _0(a_0 \otimes a_1 \otimes ...\otimes a_n)&=&\gamma (a_0 a_1 \otimes ...\otimes a_n) \\
&=&\sum Tr(a_0 a_1 a_2^{(0)} ...\; a_n^{(0)})(a_2^{(1)}\otimes a_3^{(1)}\otimes ...\otimes a_n^{(1)})  \\
&=&\sum Tr(a_0 a_1^{(0)} \epsilon (a_1^{(1)})a_2^{(0)} ...\; a_n^{(0)})(a_2^{(1)}\otimes a_3^{(1)}\otimes ...\otimes
 a_n^{(1)})    \\
&=&\sum Tr(a_0 a_1^{(0)} a_2^{(0)} ...\; a_n^{(0)})(\epsilon (a_1^{(1)})a_2^{(1)}\otimes a_3^{(1)}\otimes ...\otimes
 a_n^{(1)})\\
&=&\delta  _0 \circ \gamma (a_0 \otimes a_1 \otimes ...\otimes a_n). 
\end{eqnarray*}
For $1\leq i \leq n-1,$ it is obvious that $\gamma$  commutes with $\delta _i.$ For $i=n $ we  have \\
\begin{eqnarray*}
\gamma \circ\delta  _n(a_0 \otimes a_1 \otimes ...\otimes a_n)&=&\gamma (a_na_0 \otimes a_1
\otimes  ...\otimes a_{n-1})\\
&=&\sum Tr(a_na_0 a_1^{(0)} a_2^{(0)} ...\; a_{n-1}^{(0)})(a_1^{(1)}\otimes a_2^{(1)}\otimes ...\otimes
 a_{n-1}^{(1)})\\
&=&\sum Tr(a_0 a_1^{(0)} a_2^{(0)} ...\; a_n^{(0)})(\delta  (a_n^{(1)})a_1^{(1)} 
\otimes a_2^{(1)} \otimes ... \otimes
 a_{n-1}^{(1)})\\
&=&\delta  _n \circ \gamma (a_0 \otimes a_1 \otimes ...\otimes a_n),
\end{eqnarray*}
 where we have made use of the $\delta  $-trace property of $Tr$. \\
We leave it to the reader to check $\gamma \circ \sigma _i = \sigma _i \circ \gamma$.
Finally, we show that $\gamma$ commutes with $\tau _n$: \\
\begin{eqnarray*}
\gamma \circ \tau _n (a_0 \otimes a_1 \otimes ...\otimes a_n) &=& \gamma 
(a_n \otimes a_0 \otimes  ...\otimes a_{n-1}) \\
&=&\sum Tr (a_na_0^{(0)} a_1^{(0)} a_2^{(0)} ...\; a_{n-1}^{(0)})(a_0^{(1)}
\otimes a_1^{(1)}\otimes ...\otimes
 a_{n-1}^{(1)})
\end{eqnarray*}
and
\begin{eqnarray*}
 \tau _n \circ \gamma (a_0 \otimes a_1 \otimes ...\otimes a_n)
&=&  \sum Tr (a_0 a_1^{(0)} a_2^{(0)}  ...\; a_n^{(0)}) \tau 
(a_1^{(1)}\otimes a_2^{(1)}\otimes ...\otimes a_n^{(1)})\\
&=& \sum Tr (a_0a_1^{(0)} a_2^{(0)}  ...\; a_n^{(0)})
 ( S_{\sigma} ( a_1^{(1)} a_2^{(1)} ...a_n^{(1)} )\\ && \otimes a_1^{(2)} \otimes ...\otimes\delta 
  (a_n^{(2)})a_{n-1}^{(2)})\\
&=&\sum Tr (a_0^{(0)}a_1^{(0)} a_2^{(0)}  ...\; a_n^{(0)})\;  (a_0^{(1)}\otimes a_1^{(1)}\otimes ...  \otimes\delta 
 (a_n^{(1)})a_{(n-1)}^{(1)})\\
&=&\sum Tr (a_na_0^{(0)}a_1^{(0)} a_2^{(0)}  ...\; a_{(n-1)}^{(0)})(a_0^{(1)}\otimes a_1^{(1)}\otimes ...\otimes
 a_{(n-1)}^{(1)}).
\end{eqnarray*}
\end{proof}
\begin{corollary}
Under the  conditions of Proposition  2.1,$ \gamma $ induces the following canonical map: \\
\begin{equation} 
\gamma_\ast :\; HC_\ast (A)\longrightarrow {\widetilde{HC}_\ast}^{(\delta  ,\sigma)}(\mathcal{H}).
\end{equation}
\end{corollary}
\begin{example}
Let $A$ be $sl_q(2 , \mathbb{C})$ and $\mathcal{H }$ be $\mathbb{C} [ z , z^{-1}].$
Then $\mathcal{H }$ has a natural coaction on $A$ as follows. 
If we denote the generators of $sl_q(2 , \mathbb{C})$ by $a , b , c , d, $ then our coaction is 
\begin{eqnarray*}
\beta (a)&=& a \otimes z \\
\beta (b)&=& b \otimes z^{-1}\\ 
\beta (c)&=& c \otimes z \\
\beta (d)&=& d \otimes z ^{-1}.
\end{eqnarray*} 
If we consider
\begin{eqnarray*} 
Tr(x)=\left \{ \begin{array}{ll}
1 & \textrm{if $ x=a^kd^k $ for $k \geq 0 $} \\
0 & \textrm{otherwise}
\end{array}\right.
\end{eqnarray*} \\
then $Tr$ is a $ \epsilon$-trace and it is   $1$-invariant.
\end{example}

\begin{prop}
Let  $\mathcal{H}$ be a Hopf algebra with $S_{\sigma}^2=id$ for some grouplike element $\sigma$. Then the following
map is a cyclic map:
$$\theta : \widetilde{\mathcal{H}}^{(\epsilon ,\sigma )}_n \longrightarrow C_n (\mathcal{H}),$$
$$\theta (h_1 \otimes  h_2 \otimes  \dots \otimes  h_n ) = (S_\sigma (h_1^{(1)} h_2^{(1)} \cdots h_n^{(1)} )\otimes
 h_1^{(2)}\otimes  
h_2^{(2)} \otimes  \dots \otimes  h_n^{(2)} ).$$
Here, $C_n (\mathcal{H})$ is the corresponding cyclic module when $\mathcal{H}$ is considered 
to be an
algebra.
\end{prop}

\begin{proof}
We must show that $\theta$ commutes with $\delta _i ,\; \sigma_i$ and $\tau$. In this proof we show
$\theta\delta _n =\delta _n \theta ,\; \theta\tau = \tau \theta$ and leave the other cases to the reader.\\
 At first, let $h\in \mathcal{H}$,
so $h= S_\sigma^2 (h)= \sigma S^2 (h)\sigma^{-1} $.  Therefore $h^{(2)}S_\sigma (h^{(1)} ) =
\sigma S^2 (h^{(2)} )S(h^{(1)} ) = \sigma \epsilon (h)$. Now,

\begin{eqnarray*}
\theta\delta _n (h_1 \otimes  h_2 \otimes  \dots \otimes  h_n ) &=& \epsilon (h_n) \theta (h_1\otimes   h_2 \otimes 
 \dots \otimes  h_{n-1} ) \\
&=& \epsilon (h_n)(S_\sigma (h_1^{(1)}\cdots h_{n-1}^{(1)} ) \otimes  h_1^{(2)} \otimes  h_2^{(2)}\otimes  \dots
 \otimes  h_{n-1}^{
(2)} ) \\
\delta _n \theta (h_1 \otimes  h_2 \otimes  \dots \otimes  h_n ) &=&\delta _n (S_\sigma (h_1^{(1)}\dots h_n^{(1)}
 )\otimes   h_1^2 \otimes  \dots
\otimes  h_n^{(2)} ) \\
&=& (h_n^{(2)} S_\sigma (h_1^{(1)}\dots h_n^{(1)} )\otimes   h_1^{(2)} \otimes  \dots \otimes  h_{n-1}^{(2)} ) \\
&=& \epsilon (h_n) (S_\sigma (h_1^{(1)} \cdots h_{n-1}^{(1)} )\otimes  h_1^{(2)} \otimes  \dots \otimes  h_n^{(2)}),
\end{eqnarray*}
{so that  $\theta\delta _n =\delta _n\theta.$  Next, we have  } 
\begin{eqnarray*}
\theta \tau (h_1 \otimes  h_2 \otimes  \dots \otimes  h_n ) &=& \theta (S_\sigma (h_1^{(1)} \cdots h_n^{(1)} )\otimes 
 h_1^{(2)} \otimes  
\dots \otimes  \epsilon (h_n^{(2)} )h_{n-1}^{(2)} ) \\
&=& \epsilon (h_n^{(3)} )(S_\sigma (S_\sigma (h_1^{(2)} \cdots h_n^{(2)} )h_1^{(3)} \cdots
h_{n-1}^{(3)})\\&& \otimes  S_\sigma (h_1^{(1)} \cdots h_n^{(1)} )\otimes  h_1^{(4)} \otimes  \dots \otimes 
 h_{n-1}^{(4)} ) \\
&=& \epsilon (h_n^{(3)} )(S_\sigma^2 (h_n^{(2)} )\otimes  S_\sigma (h_1^{(1)} \cdots h_n^{(1)} )\otimes  
h_1^{(2)} \\ && \otimes  h_2^{(2)} \otimes  \cdots \otimes  h_{n-1}^{(2)} ) \\
&=& (h_n^{(2)} \otimes  S_\sigma (h_1^{(1)} \cdots h_n^{(1)} )\otimes   h_1^{(2)} \otimes  \cdots \otimes 
 h_{n-1}^{(2)} ).
\end{eqnarray*}
On the other hand,
\begin{eqnarray*}
\tau\theta (h_1 \otimes  h_2 \otimes  \dots \otimes  h_n) &=& \tau (S_\sigma (h_1^{(1)} \cdots h_n^{(1)} )\otimes 
 h_1^{(2)} \otimes  \dots \otimes  h_n^{(2)} ) \\
&=& (h_n^{(2)} \otimes  S_\sigma (h_1^{(1)} \cdots h_n^{(1)} )\otimes   h_1^{(2)} \otimes  \dots \otimes 
 h_{n-1}^{(2)} ).
\end{eqnarray*}
\end{proof}
One knows that any Hopf algebra  $\mathcal{H}$ has a right
coaction on itself by comultiplication. Let $\mathcal{H}$ have a $\sigma$-invariant
 trace $Tr$. By Proposition 3.1, we
 have a map
$$\gamma : C_n (\mathcal{H}) \longrightarrow \widetilde{\mathcal{H}}^{(\epsilon ,\sigma)}_n,$$
$$\gamma (h_0 \otimes  h_1 \otimes  \dots \otimes  h_n )= \sum Tr (h_0h_1^1 \cdots h_n^1)(h_1^2 \otimes  \dots \otimes 
 h_n^2 ).$$
\begin{theorem}
Let $\mathcal{H}$ be a Hopf algebra with $S_\sigma^2 = id$ for some $\sigma$, and also, 
let $\mathcal{H}$ have a $\sigma$-invariant trace, Tr, and let Tr($\sigma$)
be invertible in $k.$ Then, $\widetilde{HC}_n^{(\epsilon, \sigma)}(\mathcal{H})$ is a direct summand in  
${HC}_n(\mathcal{H})$ where $\mathcal{H}$ is considered  as  an algebra.
\end{theorem}
\begin{proof}
It can be shown that $ \gamma  \theta = Tr (\sigma )id $.
\end{proof}
\begin{example}
Let $G$ be a discrete group and $\mathcal{H} = kG$ be its group algebra over $k.$ Let $\sigma$
 be any central element in $G$.  It is obvious that  the following
 trace satisfies all conditions needed in the previous theorem: 
\begin{eqnarray*} 
Tr(x)=\left \{ \begin{array}{ll} 
1 & \textrm{if $ x=\sigma $ } \\
0 & \textrm{otherwise.}
\end{array}\right.
\end{eqnarray*}
\end{example}
\section{ Relation with Hopf Algebra Homology}
\textrm{
In this section we recall the  analogue of group homology for Hopf algebras and relate our cyclic homology, for 
 cocommutative Hopf algebras,  to this homology. Let $\mathcal{H}$ be a Hopf
 algebra and $ M $ a left $\mathcal{H}$-module. We define two new modules, the module of invariants and coinvariants,
   $ M^{\mathcal{H}} $ and
$M_{\mathcal{H}}$ by :}
\begin{equation}
M^{\mathcal{H}} = \{ m\in M | hm=\epsilon (h)m \hspace{15pt} \forall \; h \in \mathcal{H} \}
\end{equation}
\begin{equation}
M_\mathcal{H} = M/\{\textrm{submodule generated by} (hm-\epsilon (h)m) | h\in \mathcal{H},
m\in M \}. 
\end{equation}
\textrm{
In fact, we have two functors $-^{\mathcal{H}}$ and $-_{\mathcal{H}}$ from $\mathcal{H}$-mod
to $k$-mod and by recalling that a trivial $\mathcal{H}$-module is an $\mathcal{H}$ module
where $hm = \epsilon (h)m$ for all $m\in M,$ and all $h\in \mathcal{H}$, we see that  $M^\mathcal{H}$ 
is the biggest trivial submodule of $M$, and $M_{\mathcal{H}}$ is the biggest quotient module of $M$ that is
trivial under the action of $\mathcal{H}$. On the other hand, we have the trivial module functor from 
$k$-mod to $\mathcal{H}$-mod with $-_{\mathcal{H}}$ as its right adjoint and 
$-^\mathcal{H}$  its left adjoint. It is obvious that $M_\mathcal{H} = k \underset{\mathcal{H}}{\otimes} M$ and
 $M^\mathcal{H} = \hom_\mathcal{H} (k,M)$ where $k$ is the trivial $\mathcal{H}$-module. In the following $L_{\ast}$
  and $R^{\ast}$ denote the left and right derived functors.}
\begin{definition}
Let $M$ be an $\mathcal{H}$-module. We define  $H_\ast (\mathcal{H} ;M)$ to be  $L_\ast (-^\mathcal{H} )
(M)$ and call them  Hopf algebra homology groups of $\mathcal{H}$ with coefficients in $M$.
Using the above notations, we have $H_\ast(\mathcal{H} ;M) \cong Tor_\ast^\mathcal{H}(k ,M).$
 Similarly, if we define $H^\ast (\mathcal{H} ;M)$ to be $R^\ast (-_\mathcal{H} )(M),$ we have
$H^\ast (\mathcal{H} ; M) \cong Ext_\mathcal{H}^\ast (k,M).$
\end{definition}
\begin{example}
Let $\mathfrak{g}$ be a Lie algebra and $\mathcal{H}=U(\mathfrak{g})$ its enveloping algebra. Then, an
 $\mathcal{H}$-module
is exactly a $\mathfrak{g}$-module and we have $M^\mathcal{H} = M^{\mathfrak{g}}$ and  $M_\mathcal{H} = M_{{g}}.$ So
 $H^{\ast}(\mathcal{H} ; M) = H^{\ast}(\mathfrak{g} ; M)$ is the Lie algebra cohomology  and, $H_{\ast}(\mathcal{H} ;
 M) = H_{\ast}(\mathfrak{g} ; M).$ Similarly if $\mathcal{H}=kG$ is the group algebra of a (discrete) group $G$, then
  $H_{\ast}(\mathcal{H};M)\cong H_{\ast}(G;M)$ is the group homology and $H^{\ast}(\mathcal{H};M)\cong H^\ast(G;M).$ 
\end{example}
For every simplicial object $M$ one can define its path space, $EM$, where $EM_n = M_{n+1}$ and 
its $n^{th}$ face is $(n+1)^{th}$ face of $M$ and the same for degeneracies.
So if one denotes the path space of $\widetilde{\mathcal{H}}^{(\delta  , \sigma )}$  by $E\mathcal{H}$,  its
 simplicial structure is:
\begin{eqnarray*}
  {\delta }_i (h_0 \otimes h_1 \otimes 
 ... \otimes h_n) &=& h_0\otimes h_1 \otimes ...\otimes 
h_i h_{i+1}\otimes ... \otimes h_n \;\; 0 \leq i \leq n-1 \\
{\delta }_n(h_0 \otimes h_1 \otimes  ... \otimes h_n) &=&\delta  (h_n)h_0 \otimes h_1 \otimes ... \otimes h_{n-1} \\
 {\sigma}_i(h_0 \otimes h_1 \otimes ... \otimes h_n) &=& h_0\otimes h_1 ...\otimes h_i \otimes 1 \otimes
 h_{i+1} ... \otimes h_n  \;\;  0\leq i \leq n-1 \\
 {\sigma }_n(h_0 \otimes h_1 \otimes  ... \otimes h_n) &=& h_0 \otimes h_1  \otimes ...\otimes 1. 
\end{eqnarray*} 
It is easy to verify that $E \mathcal{H}$ is a simplicial $k$-module,  contractible and  also a free resolution for
 $k$
 via  $\delta  : E \mathcal{H} \longrightarrow  k$.

Now let $M_{\ast}$ be a chain complex of $\mathcal{H}$-modules. We denote the hyper derived  functors
$\mathbb{L}(-_{\mathcal{H}})(M_{\ast})$  by $\mathbb{H}(\mathcal{H} ; M_{\ast})$ and 
$\mathbb{R}(-^{\mathcal{H}})(M^{\ast})$  by $\mathbb{H}(\mathcal{H} ; M^{\ast}),$
where $M^{\ast}$ is a cochain complex of $\mathcal{H}$-modules.
\begin{lem}
If $\mathcal{H}$ is a cocommutative Hopf algebra then $E \mathcal{H}$ is a cyclic $k$-module 
 with $$ t_n(h_0 \otimes 
\dots \otimes h_n) = \sum (h_0 h_1^{(1)} \dots h_n^{(1)}) \otimes S(h_1^{(2)} \dots h_n^{(2)})
\otimes h_1^{(3)} \otimes  \dots \otimes h_{n-1}^{(3)}.$$
\end{lem}
\begin{proof}
As always it is needed to verify  the relations  2.2,...,2.6. We only check 2.2 and leave the others to the reader.
$$t_n^2 (h_0 \otimes 
\dots \otimes h_n) =t_n( \sum (h_0 h_1^{(1)} \dots h_n^{(1)}) \otimes S(h_1^{(2)} \dots h_n^{(2)})
\otimes h_1^{(3)} \otimes  \dots \otimes h_{n-1}^{(3)} ) $$
$$=\sum (h_0h_1^{(1)}h_2^{(1)}\dots 
h_n^{(1)}S(h_1^{(4)}h_2^{(4)}\dots 
h_n^{(4)}) h_1^{(5)}\dots \otimes h_{n-1}^{(5)}\otimes $$
$$ S(S(h_1^{(3)}h_2^{(3)} \dots h_n^{(1)})h_1^{(6)}\dots h_{n-1}^{(6)})\otimes 
S(h_1^{(2)}\dots h_n^{(2)})\otimes h_1^{(7)}\otimes \dots \otimes h_{n-2}^{(7)}$$
$$=\sum (h_0h_1^{(1)}\dots h_{n-1}
^{(1)}\otimes h_n^{(1)}\otimes
S(h_1^{(2)}h_2^{(2)}\dots h_n^{(2)} )
\otimes h_1^{(3)}\otimes \dots \otimes
h_{n-2}^{(3)}. $$
By a similar argument we get,
$$t_n^{n}(h_0 \otimes \dots  \otimes h_n) = 
\sum h_0h_1^{(1)}\otimes h_2^{(3)}\dots \otimes h_n^{(3)} 
\otimes S(h_1^{(2)}\dots h_n^{(2)}), $$
and finally we have 
$$t_n^{n+1}(h_0 \otimes \dots  \otimes h_n)= h_0 \otimes \dots  \otimes h_n. $$
\end{proof}
\textrm{
From now on we denote} $B \mathcal{H}$ for $\widetilde{\mathcal{H}}^{(\epsilon, 1 )}.$
\begin{lem} The projection $\pi : E \mathcal{H} \longrightarrow  B \mathcal{H} $
where\\
$\pi(h_o \otimes \dots \otimes h_n) = (\epsilon (h_o)h_1 \otimes \dots \otimes h_n)$
is a simplicial map and, if $\mathcal{H}$ is  cocommutative,  then, $\pi$ is a cyclic map.
\end{lem}

\begin{proof}
We leave it to the reader the first part of proof and just prove the second part.
We must verify that the following diagram is commutative.
$$
\CD
E \mathcal{H}_n @>\pi >> B \mathcal{H}_n \\
@VtVV   @VV\tau V \\
E \mathcal{H}_n @>\pi >> B \mathcal{H}_n. \\
\endCD
$$ 
We have
$$ \tau \circ \pi(h_0 \otimes \dots \otimes h_n) = \tau(\epsilon (h_0)h_1
\otimes \dots \otimes h_n) $$
$$= \epsilon (h_0)\sum (S(h_1^{(1)}h_2^{(1)}\dots h _{n-1}^{(1)}h_n )\otimes h_1^{(2)}\otimes   \dots   \otimes
 h_{n-1}^{(2)}), $$
and 
$$   \pi \circ t (h_0 \otimes \dots \otimes h_n)=\pi(\sum(h_0h_1^{(1)}\dots h_n^{1}) 
 \otimes S(h_1^{(2)}h_2^{(2)} \dots h_n^{(2)})\otimes h_1^{(3)}\otimes \dots
h_{n-1}^{(3)}$$
$$=\sum ( \epsilon (h_0h_1^{(1)}h_2^{(1)}\dots h_n^{(1)})S(h_1^{(2)}\dots h_n^{(2)})\otimes h_1^{(3)}\otimes \dots 
 h_{n-1}^{(3)} ) $$
$$ = \epsilon (h_0)\sum (S(h_1^{(1)} \dots h_{n-1}^{(1)}h_n)\otimes h_1^{(2)}\otimes \dots \otimes h_{n-1}^{(2)}). $$
\end{proof}

It is obvious that $E \mathcal{H}$ is an $\mathcal{H}$-module via  $h(h_0\otimes \dots \otimes h_n) = (hh_0\otimes
 \dots \otimes h_n)$, and the relation between $E \mathcal{H}$ and $B \mathcal{H}$ becomes
 $B \mathcal{H}=k\otimes _{\mathcal{H}}E \mathcal{H}$.\\

Next, we prove a theorem  which computes the cyclic homology of cocommutative Hopf algebras. For $\mathcal{H} = kG$,
 our result reduces to Karoubi's theorem $~\cite{kar}$.
\begin {theorem}
If $\mathcal{H}$ is a cocommutative Hopf algebra then
$$ \widetilde{HC}^{(\epsilon,1)}_n(\mathcal{H}) = \bigoplus _{i\geq 0} H_{n-2i}(\mathcal{H} ; k).$$
\end{theorem}

\begin{proof}
By the above remark we have $CC{\ast \ast}(B \mathcal{H})\cong k \otimes 
_{\mathcal{H}}CC_{\ast \ast }(E\mathcal{H}) $ where $CC_{\ast \ast} $
denotes cyclic double complex $~\cite{ld}$. Since $ E \mathcal{H} $ is contractible and 
$\epsilon : E \mathcal{H} \longrightarrow  k $ is a Homotopy equivalence the double
 complex $CC_{\ast \ast}(E \mathcal{H})$ is a resolution 
for  $ k_{\ast} $, where $ k{\ast}$ is
 $$k \leftarrow 0 \leftarrow  
 k \leftarrow 0 \dots . $$ 
On the other hand, $$ \widetilde{HC}^{(\epsilon,1)}_n(\mathcal{H})  = 
H_n(TotCC_{\ast \ast }(B \mathcal{H})) = $$
$$
 H_n(TotCC_{\ast \ast }(k\otimes _{\mathcal{H}}E\mathcal{H}) 
= \mathbb{H}(\mathcal{H} ; k_{\ast}). $$ 
So, to complete the proof it suffices to
 compute $ \mathbb{H}(\mathcal{H} ; k_{\ast}). $ 
But by finding a Cartan-Eilenberg resolution 
for $k_\ast$, we have
\[\mathbb{H}(\mathcal{H} ; k_{\ast}) = \bigoplus
 _{i\geq 0} H_{n-2i}(\mathcal{H} ; k).\]
\end{proof}

\begin{example}
Let $\mathfrak{g}$ be a Lie algebra over $k$ and $U(\mathfrak{g})$ be 
its enveloping algebra. One knows that 
$H_n(U(\mathfrak{g}) ; k) = H_n(\mathfrak{g} ; k),$ $~\cite{ld}$
so by Theorem 4.1 we have  
$$\widetilde{HC}^{(\epsilon,\sigma)}_n(U(\mathfrak{g})) = \bigoplus _{k\geq 0}H_{n-2k}(\mathfrak{g} ; k). $$
\end{example}
\textrm{
Now let $\mathcal{H}$
be a Hopf algebra  and $M$ be an $\mathcal{H}$-bimodule. We can convert $M$
to a new left $\mathcal{H}$-module, $\widetilde{M} = M$, where the action of $\mathcal{H}$ is }
$$ h \blacktriangleright m= h^{(2)}mS(h^{1}).$$

\begin{prop}{(Mac Lane isomorphism for Hopf algebras)}\\
Under the above hypotheses there is a canonical isomorphism 
$$\theta _{\ast} : H_n (\mathcal{H} , M) \cong H_n(\mathcal{H} ;
\widetilde{M}).$$
\end{prop}

\begin{proof}
If $C_n(\mathcal{H} ; \widetilde{M}) = 
\mathcal{H}^{\otimes n} \otimes \widetilde{M}$ 
then it is obvious that $C_n(\mathcal{H} ; \widetilde{M})$ is a simplicial module by the following faces and
 degeneracies:
\begin{eqnarray*}
\delta  _0 (h_1 \otimes \dots \otimes h_n\otimes m) &=& 
(\epsilon (h_1)h_2 \dots \otimes h_n \otimes m) \\
\delta  _i (h_1 \otimes \dots \otimes h_n\otimes m) &=&
(h_1 \otimes  \dots \otimes h_ih_{i+1} \otimes \dots \otimes h_n\otimes m) \hspace{1cm} 1\leq i \leq n-1 \\
\delta  _n(h_1 \otimes \dots \otimes h_n\otimes m) &=& ( h_1 \otimes \dots \otimes h_{n-1}\otimes
 h_n \blacktriangleright m) \\
 \sigma_i (h_1 \otimes  \dots \otimes h_n\otimes m ) &=& (h_1 \otimes \dots \otimes h_i\otimes 1 \otimes h_{i+1} \dots
  \otimes h_n \otimes m)  \hspace{1cm} 0 \leq i \leq n. \\
\end{eqnarray*}
Now, $H_{\ast}(\mathcal{H} ; \widetilde {M}) $ 
can be computed by the above complex. Let 
$$\theta (m \otimes h_1 \otimes \dots  \otimes h_n) = 
\sum(h_1^{(2)}\otimes  \dots \otimes h_n^{(2)}\otimes mh_1^{(1)}h_2 ^{(1)}\dots h_n^{(1)}). $$
We leave to the reader to show $\theta $ is a 
simplicial map and in fact, $\theta$ is a simplicial 
isomorphism by the following inverse map
$$\theta ^{-1}(h_1 \otimes \dots  
\otimes h_n \otimes m ) = \sum(mS(h_1^{(1)}h_2^{(1)}
\dots h_n^{(1)}) \otimes h_1^{(2)}\otimes 
\dots \otimes h_n^{(2)}). $$
\end{proof}
\textrm{
Let $k$ be the trivial $\mathcal{H}$-module.
 Then $\widetilde{k}$ is also a trivial 
 $\mathcal{H}$-module. So we have, 
$$H_n(\mathcal{H} , k) = H_n(\mathcal{H} ; k), $$
where the right hand side is the Hochschild homology 
of $\mathcal{H}$ with trivial coefficients via $\epsilon$ for
both  left and right action of $\mathcal{H}$ 
on $k$. The left hand side is the Hopf algebra
 homology of $\mathcal{H}$ via $\epsilon$.}

\begin{corollary}
If $\mathcal{H}$ is a cocommutative Hopf algebra then, 
$$\widetilde{HC}^{(\epsilon,1)}_n(\mathcal{H}) =
 \bigoplus _{k\geq 0} H_{n-2k}(\mathcal{H} , k). $$
\end{corollary}

\begin{example}
Let $V$ be a $k$-module
 and $T(V)$ be its tensor
 algebra .Then $T(V)$ is a cocommutative Hopf algebra. We 
have $H_0(T(V),k)=k,  H_1(T(V),k)=V  $ and the other
 homology groups are zero so,
$$ \widetilde{HC}^{(\epsilon,1)}_n (T(V)) = k \hspace {2cm} \textrm{if $n$ is even} $$
 $$\hspace{-2mm} \widetilde{HC}^{(\epsilon,1)}_n (T(V)) = V \hspace {2cm} \textrm{if $n$ is odd.} $$ 
\end{example}
{\bf Remark} (Connes-Moscovici cyclic cohomology of commutative Hopf algebras).\\
Using methods similar to the  above, one can compute  Connes-Moscovici periodic 
cyclic cohomology $HP^\ast_{(\epsilon,1)}(\mathcal{H})$ of commutative Hopf 
algebras. Since proofs are similar we only indicate the main steps.  
\begin{prop}
Let $\mathcal{H}$ be a commutative Hopf algebra. 
Let $(E\mathcal{H})_n=\mathcal{H}^{\otimes{n+1}},\;n\geq 0$. Then the following operators 
define a cocyclic module structure on $E \mathcal{\mathcal{H}}$: 
\begin{eqnarray*}
d_i(h_0 \otimes  \dots \otimes  h_n)&=& h_0 \otimes  \dots \otimes \Delta(h_i)\otimes \dots \otimes h_n\\
d_{n+1}(h_0 \otimes \dots\otimes h_n)&=& h_0 \otimes \dots \otimes h_n \otimes 1\\
s_i(h_0 \otimes \dots\otimes h_n)&=& h_0 \otimes \dots \otimes \epsilon(h_i)\otimes \dots \otimes h_n\\
\tau(h_0 \otimes \dots\otimes h_n)&=& h_0^{(1)}\otimes h_0^{(2)}S(h_1^{(n-1)})h_2 \otimes  \dots
 \otimes h_0^{(n)}S(h_1^{(1)}) h_n.
\end{eqnarray*}
$~~~\shoveright{\square}$
\end{prop}
Let $\mathcal{H}_{(\epsilon,\delta)}$ denote the Connes-Moscovici cocyclic
 module of $\mathcal{H}$ for the modular pair $(\epsilon,1)$.
  \begin{prop}
The following map is a morphism of cocyclic modules,
 $\psi:  \mathcal{H}^{(\epsilon,1)} \longrightarrow E \mathcal{\mathcal{H}}$,\\
\begin{center}
$\psi(h_1 \otimes\dots\otimes h_n)= 1 \otimes h_1 \otimes \dots \otimes h_n.$ 
\end{center}
$~~~\shoveright{\square}$
\end{prop}
Using the above two propositions and dualizing the above method  to prove a similar result for cocommutative Hopf
 algebras, one can prove:  

\begin{theorem}
Let $\mathcal{H}$ be a commutative Hopf algebra. Its periodic cyclic cohomology in the sense of 
Connes-Moscovici~\cite{aChM99,aChM98}  
 is given by $$HP^n_{(\epsilon,1)}(\mathcal{\mathcal{H}})=\bigoplus_{i=n~(\text{mod}~ 2 )}
  H^i(\mathcal{\mathcal{\mathcal{H}}},k)	$$
$~~~\shoveright{\square}$
\end{theorem}

For example, if $\mathcal{\mathcal{H}}=k \lbrack G \rbrack$ is the algebra 
of regular functions on an affine algebraic group $G$, 
then the coalgebra complex of $\mathcal{\mathcal{H}}=k \lbrack G \rbrack$ is
 isomorphic to the group cohomology complex, 
with trivial coefficient, where 
instead of arbitrary cochains one uses regular functions $G \times G \times \dots \times G \rightarrow k $. Denote 
this cohomology by $H^i(G,k)$. 
It follows that $$HP^n_{(\epsilon,1)}(k \lbrack G \rbrack)=\bigoplus_{i=n~~(\text{mod}~2)}H^i(G;k).$$  
For $k=\mathbb{R}$, this gives  an alternative proof of Prop.$4$ and Remark $5$ in ~\cite{achm00}.

\section{The Cyclic Homology of $A(SL_q(2,k)$ and $U_q(sl(2,k))$}

In this section we compute our cyclic homology theory for $A(SL_q(2,k))$, 
the quantized algebra of functions on the quantum group $SL_q(2,k)$ and
 also for the quantized universal enveloping algebra $U_q(sl(2,k))$.

 Let $k$ be a field of characteristic zero and $q\in k$, $q\neq 0$ and $q$ 
not a root of unity. The Hopf algebra $\mathcal{H}=A(SL_q(2,k))$ is defined as follows.
 As an algebra it is generated by symbols $x,\; u,\; v,\; y,$ with the following relations:
$$ ux=qxu, \;\; vx=qxv, \;\; yu=quy,\;\;yv=qvy,$$
$$uv=vu,\;\; xy-q^{-1}uv=yx-quv=1. $$
The coproduct, counit and antipode of   $\mathcal{H}$ are defined by  
$$\Delta (x)=x \otimes x+u \otimes v,\;\;\;\Delta (u)=x \otimes u+u \otimes y, $$ 
$$\Delta (v)=v \otimes x+y \otimes v,\;\;\;\Delta (y)=v \otimes u+y \otimes y, $$
$$\epsilon (x)=\epsilon (y)=1,\;\;\;\epsilon (u)=\epsilon (v)=0,$$
$$S(x)=y,\;\; S(y)=x,\;\;S(u)=-qu,\;\;S(v)=-q^{-1}v. $$
  For more details about $\mathcal{H}$ we refer to ~\cite{kl}. Because $S^2\neq id $, to define 
our cyclic structure we  need a modular pair $(\sigma, \delta)$ in involution. Let $\delta $ be as follows:
$$\delta (x)=q,\; \; \delta(u)=0,\; \; \delta (v)=0,\; \; \delta (y)=q^{-1}.$$
And $\sigma =1$. Then we have $\widetilde {S}^2=id$.\\

For Computing cyclic homology we should at first compute the Hochschild homology $H_{\ast}(\mathcal{H}, k)$ where $k$
 is a $\mathcal{H}$-bimodule via   $\delta$, $\epsilon$ for left and right action of $\mathcal{H}$ respectively.\\
One knows $H_{\ast}(\mathcal{H}, k)=Tor_{\ast}^{\mathcal{H}^e}(\mathcal{H}, k)$, where $\mathcal{H}^e =
 \mathcal{H}\otimes \mathcal{H}^{op}$. So we need 
a resolution for $k$, or $\mathcal{H}$ as $\mathcal{H}^e$-module. We take advantage of the  free resolution for
 $\mathcal{H}$ in ~\cite{ma}. 
 
 The explicit free resolution of $A$ as $B=\mathcal{H}\otimes \mathcal{H}^o$-module is:
$$\dots \rightarrow M_2 \overset{d_2}{\rightarrow} M_1 \overset{d_1}{\rightarrow} 
M_0 \overset{\mu}{\rightarrow} \mathcal{H},$$   
  where $\mu$ is the augmentation and $M_\ast$, $\ast\ge 0$ is a family of left $B$-module 
  with their free rank over $B$ given by 
  \begin{eqnarray*}
 && rank(M_0)=1\\
 && rank(M_1)=4\\
 && rank(M_2)=7\\
 && rank(M_\ast)=8,\; \ast\ge 3.
   \end{eqnarray*}
   We give the $B$-linear differential  mapping 
   $d_\ast: M_\ast\longrightarrow M_{\ast-1}$, $\ast\ge 0$  in terms of their basis 
   over $B$. We next give the formulas whose $B$-linear 
   extensions determine the differential $d_\ast$, $\ast>0$ together with description of the 
   $B$-basis at the same time:\\
   $d_1:M_1\longrightarrow M_0=B$ \\
   is given by  
   \begin{eqnarray*}
   &&d_1(1\otimes 1\otimes e_v)=v\otimes 1-1\otimes v,\\
      &&d_1(1\otimes 1\otimes e_u)=u\otimes 1-1\otimes u,\\
	  &&d_1(1\otimes 1\otimes e_x)=x\otimes 1-1\otimes x,\\
   &&d_1(1\otimes 1\otimes e_y)=y\otimes 1-1\otimes y,
   \end{eqnarray*}
   where $e_x$, $e_y$, $e_u$, $e_v$ form a $B$-basis for $M_1$.\\
$d_2: M_2\longrightarrow M_1$
\begin{eqnarray*}
d_2(1\otimes 1\otimes (e_v\land a_u))&=&(v\otimes 1-1\otimes v)\otimes e_u-(u\otimes 1-1\otimes u)\otimes e_v,\\
d_2(1\otimes 1\otimes (e_v\land a_x))&=&(v\otimes 1-1\otimes qv)\otimes e_x-(qx\otimes 1-1\otimes x)\otimes e_v,\\
d_2(1\otimes 1\otimes (e_v\land a_y))&=&(qv\otimes 1-1\otimes v)\otimes e_y-(y\otimes 1-1\otimes qy)\otimes e_v,\\
d_2(1\otimes 1\otimes (e_u\land a_x))&=&(u\otimes 1-1\otimes qu)\otimes e_x-(qx\otimes 1-1\otimes x)\otimes e_u,\\
d_2(1\otimes 1\otimes (e_u\land a_y))&=&(qu\otimes 1-1\otimes u)\otimes e_y-(y\otimes 1-1\otimes qy)\otimes e_u,\\
d_2(1\otimes 1\otimes \vartheta_S^{(1)})&=&y\otimes 1\otimes e_x+1\otimes x\otimes 
e_y-qu\otimes 1\otimes e_v-1\otimes qv\otimes e_u,\\ 
d_2(1\otimes 1\otimes \vartheta_T^{(1)})&=& 1\otimes y\otimes e_x+x\otimes 1\otimes 
e_y-q^{-1}u\otimes 1\otimes e_v-1\otimes q^{-1}v\otimes e_u,
\end{eqnarray*}
where $\vartheta_S^{(1)}$, $\vartheta_T^{(1)}$, $e_u\land e_x$, $e_v\land e_x$,
 $e_v\land e_y $, and $e_v\land e_u$, form a $B$-basis in $M_2$.\\
 $d_3: M_3\longrightarrow M_2$\\
 is given by 
\begin{multline*}
d_3(1\otimes 1\otimes (e_v\land e_u \land e_x)) = (v\otimes 1 -1\otimes qv)
\otimes(e_u\land e_x)-\\
\shoveright{-(u\ot 1 -1 \ot qu)\ot(e_v\land e_x)
+(q^2x\ot 1-1\ot x)\ot(e_v\land e_u),}\\
~~\\
\shoveleft{d_3(1\otimes 1\otimes (e_v\land e_u \land e_y)) = (qv\otimes 1 -1\otimes v)
\otimes(e_u\land e_y)-}
\\
\shoveright{-(qu\ot 1 -1\ot u)\ot(e_v\land e_y)+(y\ot 1-1\ot q^2y)\ot(e_v\land e_u),}\\
~~\\
\shoveleft{d_3(1\otimes 1\otimes (e_v\land \vartheta _S^{(1)})) = (v\otimes 1 -1\otimes v)
\otimes \vartheta_S^{(1)}-
q^{-1}y\ot 1 \ot(e_v\land e_x)-}\\
\shoveright{-1\ot q^{-1}x\ot(e_v\land e_y)+1\ot qv\ot(e_v\land e_u),}\\
~~\\
\shoveleft{d_3(1\otimes 1\otimes (e_v\land \vartheta _T^{(1)})) = (v\otimes 1 -1\otimes v)
\otimes \vartheta_T^{(1)}-
1\ot y \ot(e_v\land e_x)-}\\
\shoveright{-x\ot 1\ot(e_v\land e_y)+1\ot q^{-1}v\ot(e_v\land e_u),}\\
~~\\
\shoveleft{d_3(1\otimes 1\otimes (e_u\land \vartheta _S^{(1)})) = (u\otimes 1 -1\otimes u)
\otimes \vartheta_S^{(1)}-
q^{-1}y\ot 1 \ot(e_u\land e_x)-}\\
\shoveright{-1\ot q^{-1}x\ot(e_u\land e_y)-qu\ot 1\ot(e_v\land e_u),}\\
~~\\
\shoveleft{d_3(1\otimes 1\otimes (e_u\land \vartheta _T^{(1)})) = (u\otimes 1 -1\otimes u)
\otimes \vartheta_T^{(1)}-
1\ot y \ot(e_u\land e_x)-}\\
{-x\ot 1\ot(e_u\land e_y)-q^{-1}u\ot 1\ot(e_v\land e_u),}
\end{multline*}
\begin{multline*}
{d_3(1\otimes 1\otimes (e_x\land \vartheta _S^{(1)})) = x\otimes 1 \otimes \vartheta_S^{(1)}-
1\ot x \otimes \vartheta_T^{(1)}-}\\
\shoveright{ - q^{-1}u\ot 1 \ot (e_v\land e_x)-1\ot v\ot(e_u\land e_x)}\\
~~\\
\shoveleft{d_3(1\otimes 1\otimes (e_y\land \vartheta _T^{(1)})) = y\otimes 1 \otimes \vartheta_T^{(1)}
-1\otimes y\otimes \vartheta_S^{(1)}-}\\
-u\ot 1\ot(e_v\land e_y)-1\ot q^{-1}v\ot(e_u\land e_y),
\end{multline*}
where $e_x\land \vartheta_S^{(1)}$, $e_y\land \vartheta_T^{(1)}$, $e_u\land \vartheta_S^{(1)}$,
 $e_u\land \vartheta_T^{(1)}$, $e_v\land \vartheta_S^{(1)}$, 
$e_v\land \vartheta_T^{(1)}$,
$e_v\land e_u\land e_x$, and $e_v\land e_u\land e_y$ form a $B$- basis in $M_3.$\\
\\
$d_{2p+4}: M_{2p+4}\longrightarrow M_{2p+3}$, $p\ge 0$\\
is given by \\
\begin{multline*}
d_{2p+4}(1\ot 1\ot (e_v\land e_u\land \vartheta_S^{(p+1)}))=(v\ot 1 -1\ot v )\ot (e_u\land \vartheta_S^{(p+1)})
-\\
\shoveright{-(u\ot 1-1\ot u)\ot (e_v\land \vartheta_S^{(p+1)})+q^{-2}y\ot 1\ot (e_v\land e_u\land e_x \land \vartheta_S^{(p)})
+
1\ot q^{-2}x \ot(e_v\land e_u\land e_y\land \vartheta_T^{(p)}),}\\
~~\\
\shoveleft{d_{2p+4}(1\ot 1\ot (e_v\land e_u\land \vartheta_T^{(p+1)}))=(v\ot 1 -1\ot v )
\ot (e_u\land \vartheta_T^{(p+1)})
-}\\
\shoveright{-(u\ot 1-1\ot u)\ot (e_v\land \vartheta_T^{(p+1)})+1\ot y 
\ot (e_v\land e_u\land e_x \land \vartheta_S^{(p)})+
x\ot 1 (e_v\land e_u\land e_y\land \vartheta_T^{(p)}),}\\
~~\\
\shoveleft{d_{2p+4}(1\ot 1\ot (e_v\land e_x\land \vartheta_S^{(p+1)}))=(v\ot 1 -1\ot qv )\ot (e_x\land \vartheta_S^{(p+1)})
-}\\
\shoveright{-qx\ot1\ot (e_v\land \vartheta_S^{(p+1)})+1\ot x \ot (e_v\land \vartheta_T^{(p+1)})
+
1\ot v\ot (e_v\land e_u\land e_x\land \vartheta_S^{(p)}),}\\
~~\\
\shoveleft{d_{2p+4}(1\ot 1\ot (e_v\land e_y\land \vartheta_T^{(p+1)}))=(qv\ot 1 -1\ot v )\ot 
(e_y\land \vartheta_T^{(p+1)})
-}\\
\shoveright{-y\ot 1\ot (e_v\land \vartheta_T^{(p+1)})+1\ot qy \ot (e_v\land \vartheta_S^{(p+1)})
+
1\ot q^{-1}v\ot (e_v\land e_u\land e_y\land \vartheta_T^{(p)}),}\\
~~\\
\shoveleft{d_{2p+4}(1\ot 1\ot (e_u\land e_x\land \vartheta_S^{(p+1)}))=(u\ot 1 -1\ot qu )
\ot (e_x\land \vartheta_S^{(p+1)})
-}\\
\shoveright{-qx\ot 1\ot (e_u\land \vartheta_S^{(p+1)})+1\ot x \ot (e_u\land \vartheta_T^{(p+1)})
-
q^{-1}u\ot 1\ot (e_v\land e_u\land e_x\land \vartheta_S^{(p)}),}\\
~~\\
\shoveleft{d_{2p+4}(1\ot 1\ot (e_u\land e_y\land \vartheta_T^{(p+1)}))=(qu\ot 1 -1\ot u )\ot 
(e_y\land \vartheta_T^{(p+1)})
-}\\
{-y\ot 1\ot (e_u\land \vartheta_T^{(p+1)})+1\ot qy \ot (e_u\land \vartheta_S^{(p+1)})
-
u\ot 1\ot (e_v\land e_u\land e_y\land \vartheta_T^{(p)}),}
\end{multline*}
\begin{multline*}
{d_{2p+4}(1\ot 1\otimes\vartheta_S^{(p+2)})=y\ot 1\ot(e_x\land\vartheta_S^{(p+1)})+1\ot 
x\ot(e_y\land \vartheta_T^{(p+1)})-}\\
\shoveright{-qu\ot 1 \ot(e_v\land\vartheta_S^{(p+1)})-1\ot qv\ot(e_u\land\vartheta_T^{(p+1)}),}\\
~~\\
\shoveleft{d_{2p+4}(1\ot 1\ot\vartheta_T^{(p+2)})=1\ot y\ot(e_x\land\vartheta_S^{(p+1)})+x\ot 
1\ot(e_y\land \vartheta_T^{(p+1)})-}\\
{-q^{-1}u\ot 1 \ot(e_v\land\vartheta_T^{(p+1)})-1\ot q^{-1}v\ot(e_u\land\vartheta_T^{(p+1)}),}
\end{multline*}
where $\omega\land\vartheta_\ast^{(0)}$ is identified with $\omega$ for 
$$\omega =e_v\land e_u\land e_x, \; \ast=S\;\;\text{or}\;\omega=e_v\land e_u\land e_y,\; \ast=T,$$
 and where $\vartheta_T^{(p+2)}$, $\vartheta_S^{(p+2)}$, $e_u\land e_x\land\vartheta_S^{(p+1)}$,
  $e_u\land e_y\land\vartheta_T^{(p+1)}$, $e_v\land e_x\land\vartheta_S^{(p+1)}$, $e_v\land 
  e_y\land\vartheta_T^{(p+1)}$,
   $e_v\land e_u\land\vartheta_T^{(p+1)}$, and $e_v\land e_u\land\vartheta_S^{(p+1)}$ 
   form a $B$-basis in $M_{2p+4}.$\\
   ~~\\
   $d_{2p+3}: M_{2p+3}\longrightarrow M_{2p+2},$  $p>0$\\
   is given by
   \begin{multline*}
   d_{2p+3}(1\ot 1\ot (e_v\land e_u\land e_x\land\vartheta_S^{(p)}))=
   (v\ot 1 -1\ot qv)\ot(e_v\land e_x\land\vartheta_S^{(p)})-\\
   \shoveright{-(u\ot 1-1\ot qu)\ot(e_v\land e_x\land\vartheta_S^{(p)})
   +q^2x\ot 1\ot(e_v\land e_u\land\vartheta_S^{(p)})-1\ot x\ot(e_v\land e_u\land \vartheta_T^{(p)}),}\\
   ~~\\
   \shoveleft{ d_{2p+3}(1\ot 1\ot (e_v\land e_u\land e_y\land\vartheta_T^{(p)}))=
   (qv\ot 1 -1\ot v)\ot(e_u\land e_y\land\vartheta_T^{(p)})-}\\
   \shoveright{-(qu\ot 1-1\ot u)\ot(e_v\land e_y\land\vartheta_T^{(p)})
   +y\ot 1\ot(e_v\land e_u\land\vartheta_T^{(p)})-1\ot q^2y\ot(e_v\land e_u\land \vartheta_S^{(p)}),}\\
   ~~\\
   \shoveleft{d_{2p+3}(1\ot 1\ot (e_v\land\vs{S}{p+1}))=(v\ot 1 -1\ot v)\ot\vs{S}{p+1}-
   q^{-1}y\ot 1 \ot (e_v\land e_x\land\vs{S}{p})-}\\
   \shoveright{-1\ot q^{-1}x\ot(e_v\land e_y\land\vs{T}{p})+1\ot qv\ot(e_v\land e_u\land\vs{S}{p}),}\\
   ~~\\
   \shoveleft{d_{2p+3}(1\ot 1\ot (e_v\land\vs{T}{p+1}))=(v\ot 1 -1\ot v)\ot\vs{T}{p+1}-
   1\ot y \ot (e_v\land e_x\land\vs{S}{p})-}\\
   \shoveright{-x\ot 1\ot(e_v\land e_y\land\vs{T}{p})+1\ot q^{-1}v\ot(e_v\land e_u\land\vs{T}{p}),}\\
   ~~\\
   \shoveleft{d_{2p+3}(1\ot 1\ot (e_u\land\vs{S}{p+1}))=(u\ot 1 -1\ot u)\ot\vs{S}{p+1}-
   q^{-1}y\ot 1 \ot (e_u\land e_x\land\vs{S}{p})-}\\
   {-1\ot q^{-1}x\ot(e_u\land e_y\land\vs{T}{p})-qu\ot 1\ot(e_v\land e_u\land\vs{S}{p}),}
   \end{multline*}
\begin{multline*}
   d_{2p+3}({1\ot 1\ot (e_u\land\vs{T}{p+1}))=(u\ot 1 -1\ot u)\ot\vs{T}{p+1}-
   1\ot y \ot (e_u\land e_x\land\vs{S}{p})-}\\
   \shoveright{-x\ot 1\ot(e_u\land e_y\land\vs{T}{p})+ q^{-1}u\ot 1\ot(e_v\land e_u\land\vs{T}{p}),}\\
   ~~\\
   \shoveleft{d_{2p+3}(1\ot 1\ot(e_x\land\vs{S}{p+1}))=x\ot 1\ot\vs{S}{p+1}
   -1\ot x\ot\vs{T}{p+1}-}\\
   \shoveright{-q^{-1}u\ot 1\ot(e_v\land e_x\land \vs{S}{p})-1\ot v\ot(e_u\land e_x\land\vs{S}{p}),}\\
   ~~\\
   \shoveleft{d_{2p+3}(1\ot 1\ot(e_y\land\vs{T}{p+1}))=y\ot 1\ot\vs{T}{p+1}
   -1\ot y\ot\vs{S}{p+1}-}\\
   {-u\ot 1\ot(e_v\land e_y\land \vs{T}{p})-1\ot q^{-1}v\ot(e_u\land e_y\land\vs{T}{p}),}
   \end{multline*}
   where $e_x\land\vs{S}{p+1}$, $e_y\land\vs{T}{p+1}$, $e_u\land\vs{T}{p+1}$, 
   $e_u\land\vs{S}{p+1}$, $e_v\land\vs{T}{p+1}$, $e_v\land\vs{S}{p+1}$, $e_v\land e_u\land e_x\land\vs{S}{p}$, and 
   $e_v\land e_u\land e_y\land\vs{T}{p}$ form a $B$-basis in $M_{2p+3}$.\\

 By a lengthy computation  $H_0(\mathcal{H},k)=0$,
  $H_1(\mathcal{H},k)= k\lbrack 1\ot e_v\rbrack \oplus k\lbrack 1\otimes e_u\rbrack$, 
   $ H_2(\mathcal{H},k)=k\lbrack 1\ot e_v\land e_x\rbrack \oplus k\lbrack 1\ot e_u\land e_x \rbrack$,
    and $H_n(\mathcal{H},k)=0$ for all $n \geq 3$. Moreover we find
   that  the operator $B=(1-\tau)\sigma N :H_1(\mathcal{H},k)\longrightarrow H_2(\mathcal{H},k)  $ 
   is bijective and we
 obtain

\begin{theorem}
For any $q\in k$ which is not a root  of unity one has\\ $\widetilde{HC}^{(\delta,1)}_1(A(SL_q(2,k)))=k\oplus k$ and
 $\widetilde{HC}^{(\delta,1)}_n(A(SL_q(2,k)))=0$ for all $n\neq 1$. \\ In particular, 
  $\widetilde{HP}^{(\delta,1)}_0(A(SL_q(2,k)))=\widetilde{HP}^{(\delta,1)}_1(A(SL_q(2,k)))=0$.  
\end{theorem}

The above theorem shows that Theorem $4.1$ is not true for non-cocommutative Hopf algebras.

The quantum universal enveloping algebra $U_q(sl(2,k))$ is an $k$-Hopf algebra which is generated as an $k$- algebra
 by symbols $\sigma$, $\sigma^{-1}$, $x$, $y$
subject to the following relations\\
$$\sigma \sigma^{-1}=\sigma^{-1} \sigma=1, \;\; \sigma x =q^2x\sigma,\;\;  \sigma y =q^{-2}y\sigma, \;\;
xy-yx = \frac{\sigma-\sigma^{-1}}{q-q^{-1}}.$$

The coproduct, counit and antipode of $U_q(sl(2,k))$ are defined by:
$$\Delta(x)=x\otimes\sigma +1\otimes x ,\;\; \Delta(y)=y\otimes 1+\sigma^{-1}\otimes y,
 \;\;\Delta(\sigma)=\sigma\otimes \sigma,$$
$$S(\sigma)=\sigma^{-1},\;\;S(x)=-x\sigma^{-1},\;\;S(y)=-\sigma
 y,$$$$\;\;\epsilon(\sigma)=1,\;\epsilon(x)=\epsilon(y)=0.$$

It is easy to check that $S^2(a)=\sigma a \sigma^{-1}$, so that  $(\sigma^{-1},\epsilon)$ is a modular pair in
 involution. As the first step to compute its cyclic homology we should find its 
Hochschild homology group with trivial coefficients. ($k$ is a $U_q(sl(2,k))$ bimodule via $\epsilon$). We define a
 free resolution for $\mathcal{H}=U_q(sl(2,k))$ as a $\mathcal{H}^e$-module as follows
\[
\begin{CD}
(*) \hspace{2cm} \mathcal{H}@<\mu<< M_0 @<d_0<< M_1 @< d_1<< M_2@< d_2<< M_3 \dots
\end{CD}
\]
Where  $M_0$ is $\mathcal{H}^e$, $M_1$ is  the free left $\mathcal{H}^e$-module generated by symbols $1\otimes 1\otimes e_\sigma,
1\otimes 1\otimes e_x,1\otimes 1\otimes e_y$, $M_2$ is the free left $\mathcal{H}^e$-module  generated by symbols $1\otimes 1\otimes e_x\land
  e_\sigma,1\otimes 1\otimes e_y\land e_\sigma,1\otimes 1\otimes e_x\land e_y $, and finally $M_3$ is generated by $1\otimes 1\otimes e_x\land
   e_y\land e_\sigma$  as a free left $\mathcal{H}^e$-module. We let $M_n=0$ for all $n\geq 4$. We claim that with the
 following boundary operators,   $(*)$ is a free resolution for $\mathcal{H}$
\begin{eqnarray*}
 &&d_0(1\otimes 1\otimes e_x)= x\otimes 1-1\otimes x \\
 &&d_0(1\otimes 1\otimes e_y)=  y\otimes 1-1\otimes y \\
 &&d_0(1\otimes 1\otimes e_\sigma)=  \sigma\otimes 1-1\otimes \sigma \\  
 &&d_1(1\otimes 1\otimes e_x\land e_\sigma)=(\sigma \otimes 1-1\otimes q^2\sigma)\otimes e_x -(q^2x\otimes1-1\otimes x
  )\otimes e_\sigma\\
&& d_1(1\otimes 1\otimes e_y\land e_\sigma)=(\sigma \otimes 1-1\otimes q^{-2}\sigma)\otimes e_y -(q^{-2}y\otimes
 1-1\otimes y )\otimes e_\sigma \\
&& d_1(1\otimes 1\otimes e_x\land e_y)=(y\otimes 1-1\otimes y)\otimes e_x-(x\otimes 1- 1\otimes x )\otimes e_y \\ 
  &&~~~~~~~~~~~~~~~~~~~  +\frac{1}{q-q^{-1}}(\sigma^{-1}\otimes \sigma^{-1}+ 1\otimes 1)\otimes e_\sigma \\
&&d_2(1\otimes 1\otimes e_x\land e_y\land e_\sigma )=(y\otimes 1-1\otimes q^2y)\otimes e_x\land e_\sigma\\
&&~~~~~~~~~ - q^2(q^2x\otimes 1-1\otimes x )\otimes e_y\land e_\sigma + q^2(\sigma\otimes 1-1\otimes \sigma)\otimes
 e_y\land e_x
\end{eqnarray*}

To show that this complex is a resolution, we need a homotopy map. First we recall that  the set
 \hbox{$\{\sigma^lx^my^n\mid l\in\mathbb{Z},m,n\in\mathbb{N}_0  \}$} is a P.B.W. type basis   for $\mathcal{H}$
  ~\cite{kl}.\\
Let $$\phi(a,b,n)=(a^{n-1}\otimes 1+a^{n-1}\otimes b+\dots + a\otimes b^{n-1}+1\otimes b^{n-1})$$  
 where $n\in \mathbb{N}, a\in \mathcal{H}, b\in \mathcal{H}^o $, and $\phi(a,b,0)=0$, and $\omega(p)=1 $ if $p \geq 0$
  and $0$ otherwise.\\
The following maps define a homotopy map  for $(*)$ i.e. $Sd+dS=1$:
\begin{eqnarray*}
&&S_{-1}: \mathcal{H} \rightarrow M_0,\\
&&S_{-1}(a)=1\otimes a, \\
&& S_0 : M_0 \rightarrow M_1,\\
&&S_0(\sigma^lx^my^n\otimes b)= (1\otimes b)((\sigma^lx^m\otimes 1)\phi(y,y,n)\otimes e_y+\\
&&~~~~~~~+(\sigma^l\otimes y^n)\phi(x,x,m)\otimes e_x)+\omega(l)(1\otimes x^my^n)\phi(\sigma,\sigma,l)\otimes e_\sigma
 \\
&&~~~~~~~~~~~~~~~~~+(\omega(l)-1)(1\otimes x^my^n)\phi(\sigma^{-1},\sigma^{-1},-l)(\sigma^{-1}\otimes
 \sigma^{-1}\otimes e_\sigma),\\
&& S_1:M_1\rightarrow M_2,\\
&& S_1(\sigma^lx^my^n\otimes b \otimes e_y)=0,\\
&& S_1(\sigma^lx^my^n\otimes b \otimes e_x)=(1\otimes b)((\sigma^lx^m\otimes 1)\phi(y,y,n)\otimes e_x\land e_y\\ 
&&+\frac{1-q^{2n}}{(q-q^{-1})(1-q^2)}(\sigma^l\otimes y^{n-1})\phi(x,x,m)(\sigma^{-1}\otimes \sigma^{-1}+q^{-2}\otimes
 1 )\otimes e_x\land e_\sigma\\
  &&~~~~~~~~+\frac{1}{q-q^{-1}}(\sigma^lx^m\otimes 1)\phi(y,y,n-1)(\sigma^{-1}\otimes\sigma^{-1}+q^2\otimes 1)\otimes
   e_y\land e_\sigma),\\
&& S_1(\sigma^lx^my^n\otimes b \otimes e_\sigma)=(1\otimes b)(q^2(\sigma^lx^m\otimes 1)\phi(y,q^2y,n)\otimes e_y\land
 e_\sigma\\
&&~~~~~~~~~~~~~~~~~~~~~~~~~~~~~~~~~~~~~~+q^{2(n-1)}(\sigma^l\otimes y^n)\phi(x,q^{-2}x,m)\otimes e_x\land e_\sigma),\\
&& S_2:M_2\rightarrow M_3,\\
&&S_2(a \otimes b \otimes e_x \land e_y)=0,\\
&&S_2(a \otimes b \otimes e_y\land e_\sigma )=0,\\
&&S_2(\sigma^lx^my^n \otimes b \otimes e_x \land e_\sigma)=(1 \otimes b)(\sigma^l x^m \otimes 1)\phi(y,q^2y,n)\otimes
 e_x \land e_y \land e_\sigma,\\
&&S_n=0:M_n \rightarrow M_{n+1}~~\text{for}~~ n\geq 3.
\end{eqnarray*}

Again, by a rather long, but straightforward computation, we can check that $dS+Sd=1$.    
By using the  definition of Hochschild homology as $Tor^{\mathcal{H}^e}(\mathcal{H},k)$ we have the following theorem
\begin{theorem}
$H_n(U_q(sl(2,k)),k)=k$ if $n=0,3$, generated by $1$ and $1\otimes e_x\land e_y\land e_\sigma$ respectively, 
 and $H_n(U_q(sl(2,k)),k)=0$ for $n\neq 0,3$. Here $k$ is a 
$U_q(sl(2,k))$-bimodule via
 $\epsilon$ for both sides. 
\end{theorem}
\begin{corollary}
$\widetilde{HC}^{(\epsilon,\sigma)}_n(U_q(sl(2,k)))=k$ when $n\neq 1$, and  $0$ for $n=1$. 
\end{corollary}


\begin{thebibliography}{9} 


   \bibitem{aChM99}

   A. Connes and H. Moscovici, \emph{Cyclic cohomology and Hopf algebras},  Moshe Flato (1937--1998). Lett. Math. Phys.
    48 (1999), no. 1, 97--108.
\bibitem{aChM98}
A. Connes and H. Moscovici, \emph{Hopf algebras, cyclic cohomology and the transverse index theorem}. Comm. Math.
 Phys. 198 (1998), no. 1, 199--246
 \bibitem{achm00}
A. Connes and H. Moscovici, \emph{ Cyclic cohomology and Hopf algebra symmetry}. Conference Moshé Flato 1999 (Dijon). Lett.
Math. Phys. 52 (2000), no. 1, 1--28. 
 \bibitem{aC85}
    A. Connes, \emph{ Noncommutative differential geometry}, Inst. Hautes 
     \`Etudes Sci. Publ. Math. No. 62, 
    257-360, (1985). 
\bibitem{cr}
M. Crainic, \emph{Cyclic cohomology of Hopf algebras}, 
J. Pure Appl. Algebra 166 (2002), no. 1-2, 29--66.
\bibitem{kar}
  M. Karoubi,\;\emph{Homologie cyclique des groups et des algebres}, C. R. Acad. Sci. Paris Ser. A-B 297(1983),381-384.
   MR 85g:18012.

\bibitem{kl}
 A. Klimyk, and  K.  Schm\"udgen, \emph{ Quantum groups and their representations}, Texts and Monographs in Physics.
  Springer-Verlag, Berlin, (1997).
\bibitem{ma}
T. Masuda, Y. Nakagami and  J. Watanabe, \emph{ Noncommutative differential geometry on the quantum ${\rm SU}(2)$. I.
 An algebraic viewpoint}, $K$-Theory 4 (1990), no. 2, 157--180.
\bibitem{ld}
     J. L. Loday, \emph{Cyclic Homology,} Springer-Verlag, (1992). 

\bibitem{sw}
   M. E. Sweedler, \emph{Hopf Algebras,} Benjamin, New York, (1969).

\bibitem{taill}
   R. Taillefer, \emph{Cyclic homology of Hopf algebras},  $K$-Theory 24 (2001), no. 1, 69--85. 
\end{thebibliography}
\end{document}